%% file: main.tex
\documentclass[a4paper,oneside,10pt]{article}
\usepackage{amsfonts}
\usepackage{graphicx} % Required for inserting images
\usepackage{mathtools}
\usepackage{relsize}
\usepackage{tikz}
\usepackage{authblk}
\usepackage{comment}
\usepackage{dirtytalk}
\usepackage{multirow}

\def\PN{{\Pi^{\nabla, E}_{k}}}
\def\PZ{{\Pi^{0, E}_{k}}}
\def\Pf{{\Pi^{0, f}_{k-1}}}
\def\PV{{\boldsymbol{\Pi}^{0,E}_{k-1}}}

\newcommand{\StabE}{\mathcal{S}^E}

\newcommand{\Hmat}{\mathbf{H}}

\newcommand{\R}{\mathbb{R}}
\newcommand{\Pk}{\mathbb{P}}
\newcommand{\Faces}{\Sigma_h}
\newcommand{\FacesE}{\Sigma^E_h}

\newcommand{\mm}{\boldsymbol{m}}
\newcommand{\nn}{\boldsymbol{n}}
\newcommand{\pp}{\boldsymbol{p}}
\newcommand{\nnE}{\nn^E}
\newcommand{\nnf}{\nn^f}
\newcommand{\xx}{\boldsymbol{x}}
\newcommand{\ii}{\boldsymbol{\imath}}
\newcommand{\jj}{\boldsymbol{\jmath}}
\newcommand{\df}{~{\rm d}f}
\newcommand{\dE}{~{\rm d}E}

\newcommand{\jump}[1]{[\![#1]\!]}

\usepackage{amsthm}
\newtheoremstyle{theorem-with-dot}% nome
  {10pt} % spazio sopra
  {10pt} % spazio sotto
  {\itshape} % corpo del teorema
  {} % indent
  {\bfseries} % intestazione (Theorem, etc.)
  {.} % PUNTO dopo l’intestazione!
  { } % spazio dopo intestazione
  {} % intestazione personalizzata
\theoremstyle{theorem-with-dot}

\newtheorem{remark}{Remark}
\newtheorem{impDet}{Implementation Detail}

\newcommand{\vmb}{\boldsymbol{\mathcal{M}}^{k-1}(E)}

\newcommand{\corr}[1]{{\color{black}#1}}

\DeclareMathOperator{\diver}{div}

\newcommand{\drawcube}[0]{ 
    \draw[dashed] (0,0,0) -- (0,2,0);
    \draw[dashed] (0,0,0) -- (2,0,0);
    \draw[dashed] (0,0,0) -- (0,0,2);
    \draw (2,0,0) -- (2,2,0);
    \draw (0,2,0) -- (2,2,0);
    \draw (0,0,2) -- (2,0,2);
    \draw (0,0,2) -- (0,2,2);
    \draw (0,2,2) -- (2,2,2);
    \draw (2,0,2) -- (2,2,2);
    \draw (2,0,0) -- (2,0,2);
    \draw (0,2,0) -- (0,2,2);
    \draw (2,2,0) -- (2,2,2);    

}

\newcommand{\redcircles}[1]{ 
    \fill[red] (1,1,0) circle (2pt);
    \fill[red] (1,1,2) circle (2pt);
    \fill[red] (1,0,1) circle (2pt);
    \fill[red] (1,2,1) circle (2pt);
    \fill[red] (0,1,1) circle (2pt);
    \fill[red] (2,1,1) circle (2pt);    
    \node[text=red, scale=0.8] at (1 - 0.15, 1 - 0.15, 0) {\small $\mathbf{#1}$};
    \node[text=red, scale=0.8] at (1 + 0.15, 1 + 0.15, 2) {\small $\mathbf{#1}$};
    \node[text=red, scale=0.8] at (1 - 0.15, 0.15, 1) {\small $\mathbf{#1}$};
    \node[text=red, scale=0.8] at (1 - 0.15, 2 + 0.15, 1) {\small $\mathbf{#1}$};
    \node[text=red, scale=0.8] at (-0.15, 1 + 0.15, 1) {\small $\mathbf{#1}$};
    \node[text=red, scale=0.8] at (2 - 0.15, 1 + 0.15, 1) {\small $\mathbf{#1}$};
}

\newcommand{\bluecircle}[1]{ 
    \fill[blue] (1,1,1) circle (2pt);
    \node[text=blue, scale=0.8] at (1 - 0.15, 1 + 0.15, 1) {\small $\mathbf{#1}$};
}

\title{The 3D Nonconforming Virtual Element Method with enhanced spaces: numerical investigation and implementation guide}
\author[1]{{{F. Dassi}} 
\thanks{franco.dassi@unimib.it}}
\author[1]{{{M. Trezzi}} 
\thanks{manuel.trezzi@unimib.it}}

\affil[1]{Dipartimento di Matematica e Applicazioni, Universit\`a degli Studi di Milano Bicocca, Via Roberto Cozzi 55 - 20125 Milano, Italy}  
\date{March 2025}

\begin{document}

\maketitle
\begin{abstract}
In this paper, we describe the 3D Nonconforming Virtual Element Method,
where the local spaces are enhanced 
in order to build an $L^2$ projection of virtual functions onto polynomials of degree $k$.
This higher order projection provides a more accurate representation of the load term and,
more importantly, it enables the construction of the reaction term.
We conduct a convergence analysis of the $L^2$ norm and $H^1$ seminorm errors 
with respect to mesh size and approximation degree.
We also compare our results with those obtained using the conforming Virtual Element Method.\newline
\phantom{ss}\newline
\noindent\textbf{Keywords:} Nonconforming Virtual Element Method; Diffusion reaction problems; Polyhedral meshes;
\end{abstract}

\input{intro}

\input{problem}

\input{vemDisc}
\input{numExe}

\smallskip
\begin{center}
{\bf Aknowledgements} 
\end{center}

\medskip
The authors have been partially funded by the European Union (ERC Synergy, NEMESIS, project number 101115663).
Views and opinions expressed are however those of the authors only and do not necessarily reflect those of the European Union or the ERC Executive Agency.

%%%%BIBLIO
\addcontentsline{toc}{section}{\refname}
\bibliographystyle{plain}
\bibliography{biblio}

\end{document}

%% file: intro.tex
% Possibili titoli:
% \begin{itemize}
%    \item A Numerical perspective on the Nonconforming Virtual Element Method for general diffusion-reaction problem in 3D.
%    \item An Investigation of the Numerical Performance of the 3D Nonconforming Virtual Element Method
%    \item The 3D Nonconforming Virtual Element Method Compared to the Conforming Approach
%    \item A Performance Analysis of the 3D Nonconforming Virtual Element Method for general diffusion-reaction problems.
%    \item Implementation and Numerical Analysis of the 3D Nonconforming Virtual Element Method
%\end{itemize}

\section{Introduction}

The Virtual Element Method (VEM) was introduced in 2013 by L.~Beir\~ao da Veiga \emph{et al.}~\cite{volley} as a generalization of the Finite Element Method. 
It allows for very general polygonal and polyhedral meshes, including non-convex shapes, and can handle regular solutions.

The key feature of the VEM lies in the definition of the basis functions. 
Unlike the FEM, where basis functions are \emph{only} polynomials of a prescribed degree,
the VEM defines basis functions as the solution of a specific Partial Differential Equation (PDE)
that includes polynomials but also other functions.
These basis functions are never computed explicitly,
and in practice, 
there is no need to know them point-wise during the matrix and right hand side assembling processes.
Furthermore, since polynomials are part of the VE space,
the accuracy of the method and the compatibility with the FEM are guaranteed.

This approach brings some advantages for both mesh discretization and the definition of functional spaces.
On one side, having the possibility to deal with generally shaped polygons/polyhedrons provides more flexibility in mesh generation.
Furthermore, if we consider mesh adaptation processes, 
the VEM allows the presence of hanging nodes, 
which facilitate the refinement process.
On the other side, the VE spaces can be built in such a way that the discrete functions inherit some useful properties.
For instance, in~\cite{vaccaDivFree},
the authors developed a family of VE spaces for Stokes problems exactly \say{divergence-free}.
FEM approaches obtain such a property in a relaxed/weak sense, 
or using special discretizations, for instance the Scott-Vogelius type elements.
Similar consideration can be made for~\cite{DASSI2020112910}
where a VEM approach was developed to solve elasticity equations based on the Hellinger–Reissner variational principle. 
The FEM imposes the symmetry of the stress tensor variable in a weak sense, 
and this constraint increases the dimension of the linear system.
The VEM discretization proposed in~\cite{DASSI2020112910} overcomes this issue by
imposing the symmetry of the tensor field within the discrete functional spaces,
resulting in a cheaper linear system, 
as no additional equations are needed to enforce the symmetry.

On of the main advantage of the nonconforming VEM over its conforming counterpart lies in its definition and implementation in the 3D case.
The conforming method requires the prior definition of a VE space on each polygonal face of a polyhedron, 
followed by the definition of the VE space in the interior.
In the nonconforming case, there is no need to define a VE space on the faces, making the overall definition more natural.
When moving from 2D to 3D, the degrees of freedom (DoFs) undergo a dimensional shift (edges become faces, and faces become volumes),
but their definition remains essentially the same, preserving the structure of the 2D setting.

B.~Ayuso de Dios \emph{et al.} introduce the nonconforming version of the VEM for second order elliptic problems~\cite{blanca}.
They provide the construction of both two and three dimensional spaces and they make the error analysis,
but they limit themselves to test these nonconforming virtual elements \emph{only} to two dimensions.
Although the work~\cite{CMS:2015} by Cangiani \emph{et al.} 
offers a valuable guide for the definition of local matrices 
in both two and three dimensions. 
To the best of our knowledge, 
the only available numerical results for the 3D nonconforming VEM are those presented in~\cite{coreani}.
Here, the authors focus on the Poisson problem and develop an analysis that accommodates polyhedra with small faces. 
Then, the nonconforming VEM has been applied to several problems~\cite{CMS:2015,LPT:2025,CMV:2016,LC:2018,ZZYC:2019,trezzi2024}.
However, all these papers provide the theory for both 2D and 3D setting, 
but they provide numerical experiments \emph{only} in 2D.

The aim of this paper is twofold: first, to numerically validate the theory established in~\cite{blanca,CMS:2015} in the three-dimensional setting; second, to provide further implementation details for the practical realization of the nonconforming VEM in 3D.

The paper is organized as follows. 
In Section~\ref{sec:model}, we present the model problem.
Then, in Section~\ref{sec:disc},
we describe the VE spaces/projections
and show how to build the linear and bilinear form for a diffusion-reaction problem.
In this section, we also outline key implementation aspects to aid and inform future developments.
Finally, in Section~\ref{sec:numExe}, we show some numerical 3D experiments.

\paragraph{Notation} 
During the paper, we adopt standard notation for Sobolev spaces and norms~\cite{adams2003pure}.

We denote by $E$ a polyhedron with
$\xx_E$, $h_E$, and $|E|$ representing its barycentre, diameter, and volume, respectively,
Similarly, we denote a generic face by $f$,
with $\xx_f$, $h_f$, and $|f|$ its barycentre, diameter and area.
Then, $\FacesE$ is the set of the faces of the element $E$ and
$\nnE$ is the outward unit normal vector to the polyhedron $E$,
similarly $\nnf$ is one of the two outward normal vectors with respect to a face $f$.

We define the jump operator for a function $v$ on a face shared by two polyhedrons as
\[
\jump{v}_f = v|_{E^+} - v|_{E^-} \, , 
\]
where $E^+$ and $E^-$ are such that $f \in \overline{E^{+}} \cap \overline{E^{-}}$
while,
if $f$ is shared by only one polyhedron, 
the jump is defined as $\jump{v}_f = v|_E$.

Then, given $\mathcal{O} \subset \mathbb{R}^d$ and a positive integer $k$,  
we denote by $\mathbb{P}_k(\mathcal{O})$ the space of polynomials in $d$ variables of degree less than or equal to $k$.  
Following common practice in the VEM,  
we choose as a basis of $\mathbb{P}_k(\mathcal{O})$ the set of scaled monomials:
$$
\mathcal{M}^s(\mathcal{O}) :=
\left\{
\left(\dfrac{\xx - \xx_{\mathcal{O}}}{h_{\mathcal{O}}} \right)^{\boldsymbol{\ell}} \, : \, |\boldsymbol{\ell}| \leq s 
\right\}\,,
$$
where $\boldsymbol{\ell}$ is a multi-index whose length depends on $d$. 
It is worth noting that when $\mathcal{O}$ is a polyhedron, 
the variable $\xx$ refers to the physical coordinates of a point in 3D space,  
whereas when $\mathcal{O}$ is a polygon (specifically, a face of a polyhedron), 
$\xx$ denotes the local coordinates on the plane defined by the polygon.

%% file: problem.tex
\section{The model problem}\label{sec:model}

Let $\Omega\subset \mathbb{R}^3$ be a polytopal domain with boundary $\Gamma$, 
we consider the following general diffusion-reaction problem with homogeneous Dirichlet boundary conditions:
\begin{equation}
\left\{
\begin{aligned}
     -\diver(\boldsymbol{\kappa}\nabla u)  + u &= f  &\text{in }\Omega \, ,\\[0.4em]
      u& = 0 &\text{on }\Gamma \, , 
    \end{aligned}
    \right.
\label{eqn:strong}
\end{equation}
where $f\in L^2(\Omega)$ is the load term, 
and $\boldsymbol{\kappa} \in [L^\infty(\Omega)]^{3\times 3}$ is the symmetric and elliptic diffusive tensor. We remark that if $\boldsymbol{\kappa}$ coincides with the identity matrix, the term $\diver(\boldsymbol{\kappa}\nabla u)$ coincides with the Laplacian operator $\Delta u$.

\begin{remark}
In the model problem, we assume that the reaction term $u$ is not multiplied by a reaction coefficient.  
This is justified by the fact that the equation can be rescaled by dividing through by the reaction coefficient, thereby simplifying the analysis without loss of generality.
\end{remark}
The variational formulation of the model problem reads as:
\begin{equation}
\left\{
\begin{aligned}
&\text{find $u\in H^1_{0}(\Omega)$ such that:} \\
&a(u,v) + m(u,v) =  F(v) \quad \forall v \in H_0^1(\Omega) \,,
\end{aligned}
\right .
\label{eqn:var}
\end{equation}
where the linear and bilinear forms are defined as
\begin{alignat}{2}
a(u, v) &:= \int_{\Omega} \boldsymbol{\kappa} \, \nabla u \cdot \nabla v \, {\rm d}\Omega 
&\quad &\text{for all $u,v \in H_0^1(\Omega)$}, \label{eq:a-c} \\
m(u, v) &:= \int_{\Omega} u \, v \,{\rm d}\Omega 
&\quad &\text{for all $u,v \in H_0^1(\Omega)$}, \label{eq:b-c}\\
F(v)    &:= \int_{\Omega} f \, v \,{\rm d}\Omega 
&\quad &\text{for all $v \in H_0^1(\Omega)$}. \label{eq:f-c}
\end{alignat}

Problem~\eqref{eqn:var} is well-posed as a consequence of the ellipticity of $\boldsymbol{\kappa}$, the coercivity of the bilinear forms $a(\cdot,\cdot)$ and $m(\cdot,\cdot)$, and the application of the Lax–Milgram lemma.

%% file: vemDisc.tex
\section{The VEM nonconforming Discretization}~\label{sec:disc}

After describing the mesh assumption, 
we introduce the nonconforming virtual element spaces and 
provide implementation details for constructing both the projection operators and the linear and bilinear forms appearing in~\eqref{eqn:var}.

\subsection{Mesh Assumptions}

Let $\{\Omega_h\}_h$ be a family of tessellations of the domain $\Omega$ into non-overlapping polyhedra.
The set of all the faces in the decomposition $\Omega_h$ is denoted by $\Faces$.
We assume that each mesh $\Omega_h$ satisfies the following properties: there exists a constant $\rho$, independent of $h$, such that

\begin{itemize}
\item each polyhedron $E$ and each face $f$ is star-shaped with respect to a ball of radius greater than or equal to $\rho h_E$ and $\rho h_f$, respectively;
\item $h_f \geq \rho h_E$ holds for every face $f\in\FacesE$.
\end{itemize}

It is worth noting that such mesh assumptions are mainly employed for the theoretical analysis of the VEM. 
In practice, however, numerical experiments have shown that the method is more robust and 
exhibits the expected behaviour even when these assumptions are not satisfied.

\begin{comment}
We introduce the following polynomial projections: %that are standard for VEM:
\begin{itemize}
\item the $\boldsymbol{L^2}$\textbf{-projection} $\Pi_n^{0, E} \colon L^2(E) \to \Pk_n(E)$, given by
\begin{equation*}
\label{eq:P0_k^E}
\int_Eq_n (v - \, {\Pi}_{n}^{0, E}  v) \, {\rm d} E = 0 \qquad  \text{for all $v \in L^2(E)$  and $q_n \in \Pk_n(E)$,} 
\end{equation*} 
with obvious extensions for  vector functions $\mathbf{\Pi}^{0,E}_{n}:[L^2(E)]^3 \to [\mathbb{P}_n(E)]^3$; or functions defined on a face $\Pi_n^{0, f}: L^2(f) \to \Pk_n(f)$;
\item the $\boldsymbol{H^1}$\textbf{-seminorm projection} ${\Pi}_{n}^{\nabla,E} \colon H^1(E) \to \Pk_n(E)$, defined by 
\begin{equation}
\label{eq:Pn_k^E}
\left\{
\begin{aligned}
& \int_E \nabla  \,q_n \cdot \nabla ( v - \, {\Pi}_{n}^{\nabla,E}   v)\, {\rm d} E = 0 \quad  \text{for all $v \in H^1(E)$ and  $q_n \in \Pk_n(E)$,} \\
& P_0(v - \PN v) = 0 \quad  \text{for all $v \in H^1(E)$} \, ,
\end{aligned}
\right.
\end{equation}
where $P_0$ is a projection onto constants, see \cite{autostoppisti} for possible choices.
\end{itemize}
\end{comment}

\subsection{The nonconforming VEM spaces}
Let $E\in \Omega_h$ be a polyhedron,
and let $k$ be a positive integer.
We introduce the enhanced nonconforming virtual element space defined as
\begin{alignat}{1}
V^{k}_h(E) 
\coloneqq
\Bigg\{ v_h \in H^1(E) : &\:  \frac{\partial v_h}{\partial \nn} \in \Pk_{k-1}(f)\quad\forall f \in \FacesE,\nonumber\\
&\: \Delta v_h \in \Pk_k(E), \nonumber\\
&\: \left(v_h - \PN v_h, \, p_k \right)_{0,E}  \! = 0\:\forall p_k \in \Pk_k(E) \setminus \Pk_{k-2}(E)\Bigg\}\, ,\nonumber\\[0.1em] 
\label{eq:local-space}
\end{alignat}
\corr{where $\Pk_k(E) \setminus \Pk_{k-2}(E)$ denotes the space of polynomials of degree exactly $k$ or $k-1$.}
Here, $\PN v_h$ is a projection operator that will be defined later in Section~\ref{sec:projNabla}. 
The dimension of this discrete space is given by the formula
\begin{equation}\label{eq:local-dimesnion}
N_E \coloneqq \dfrac{k(k+1)}{2}n_f+\dfrac{k(k-1)(k+1)}{6} \,,
\end{equation}
where $n_f$ is the number of faces of the polyhedron $E$. 
Then, to uniquely determine a function $v_h\in V^{k}_h(E)$,
we can take the following set of linear operators as Dofs:
\begin{itemize}
    \item \texttt{D1}: the face moments on each face $f\in\FacesE$:
    \begin{equation}\label{eq:dofsFaces}
    \dfrac{1}{|f|} \int_f v_h \, m_{\boldsymbol{\ell}} \, {\rm d}s \, , \quad 
    \forall m_{\boldsymbol{\ell}} \in \mathcal{M}^{k-1}(f) \, ;
    \end{equation} 
    \item \texttt{D2}: the bulk moments in $E$:
    \begin{equation}\label{eq:dofsInt}
    \dfrac{1}{|E|} \int_E v_h \, m_{\boldsymbol{\ell}} \,  {\rm d}E \, , \quad \forall m_{\boldsymbol{\ell}} \in \mathcal{M}^{k-2}(E) \, .
    \end{equation}
\end{itemize} 
\input{dofFig}
The proof that such choice of Dofs is unisolvent for $V_h^k(E)$ can be found in~\cite{blanca}.
Before showing how to use such a space to build the global one, 
we provide some comments.

First, the construction of a 3D nonconforming VE space is more natural than that of its conforming counterpart.
In the conforming case, the definition of the VE space on a polyhedron requires the introduction of a VE space also on each of its faces.
Moreover, in order to compute the projection operators, the virtual element spaces used are not those originally proposed in~\cite{volley}, 
but are enriched with the so-called enhancing property~\cite{BDR:2017}.
In contrast, as shown in Equation~\eqref{eq:local-space},
for nonconforming case there is no need to define such a space on faces:
the VE space is directly defined inside the polyhedron.

Second, a comparison with the two-dimensional counterpart~\cite{blanca} shows that 
both the local space $V_h^k(E)$ and Dofs exhibit a dimensional extension (with edges becoming faces, and faces becoming volumes) 
while preserving their definitions and the underlying structure of the 2D formulation.

Finally, the VE space introduced in Equation~\eqref{eq:local-space} differs from the one proposed in~\cite{coreani}.
The main difference lies in the regularity of the Laplacian of a generic function $v_h \in V^k_h(E)$.
In our case, we have $\Delta v_h \in \Pk_k(E)$, while in~\cite{coreani}, $\Delta v_h \in \Pk_{k-2}(E)$.
This difference stems from the presence of the enhancing property in the space \( V_h^k(E) \) defined here:
$$
\left(v_h - \PN v_h, \, p_k \right)_{E} = 0 \quad \forall p_k \in \Pk_k(E) \setminus \Pk_{k-2}(E)\,.
$$
Indeed, in order to endow the space with this property, 
its dimension must first be increased and then appropriately restricted to enforce it~\cite{BDR:2017}.

Then, the global VE space of order~$k$ is obtained by gluing the local spaces across the faces of the polyhedra:
\begin{equation}
\begin{aligned}
V^{k}_h(\Omega_h) \coloneqq
\bigl\{
v_h \in H^1(\Omega_h) \quad : \quad &
v_h|_E \in V_h^{k}(E) \quad \forall E \in \Omega_h\, , \\ 
&\int_f \jump{v_h}\, 
q \, {\rm d}f= 0 \quad  \forall q \in \mathbb P_{k-1}(f) \quad \forall f \in \Faces
\bigr\} \, .
\end{aligned}
\label{eq:global-space}
\end{equation}
This space is not globally continuous; 
instead, continuity is enforced weakly across interior faces by requiring 
that the jumps vanish when tested against polynomials of degree up to~$k-1$.

\subsection{Projectors}

A generic function $v_h \in V^k_h(E)$ is virtual, i.e., 
it is only known through its degrees of freedom. 
The VEM defines suitable projection operators 
that are essential to approximate the linear and bilinear forms in Problem~\ref{eqn:var}. 
%Following the approach of~\cite{autostoppisti, CMS:2015}, 
%for each projection operator we also provide a 
%\emph{guide} on how to compute it.
\corr{Following the approach of~\cite{autostoppisti, CMS:2015},  
this section provides a \emph{guide} on how to implement these projection operators.  
In particular, compared to the previously cited works, we place greater emphasis on the implementation details and the structure of the matrices involved.
}

\subsubsection{The projection $\Pf$}\label{sec:pfProj}

Given a face $f \in \FacesE$ of a polyhedron $E$, 
we define the projection operator $\Pf: V_h^k(E)|_f \to \mathbb{P}_{k-1}(f)$ 
as 
\begin{equation}
(p_{k-1},\, v_h - \Pf v_h)_f = 0\quad \forall p_{k-1}\in\mathbb{P}_{k-1}(f)\,,  
\label{eqn:defPZf}
\end{equation}
where $(\cdot,\cdot)_f$ denotes the $L^2$ inner product on the face $f$.

Since $\Pf v_h$ is an element of $\mathbb{P}_{k-1}(f)$, 
in Equation~\eqref{eqn:defPZf}, we can let $p_{k-1}$ vary in $\mathcal{M}^{k-1}(f)$,
and we can represent such projection operator in the basis~$\mathcal{M}^{k-1}(f)$:
$$
\Pf v_h = \sum_{|\jj|\leq k-1} q^{\jj} m^{\jj}\,.
$$
Then, Equation~\eqref{eqn:defPZf} becomes:
$$
\sum_{|\jj|\leq k-1} q^{\jj} (m^{\jj},\,m^{\ii})_f = (v_h,\,m^{\ii})_f\qquad \forall |\ii|\leq k-1\,.
$$
This linear system can be written in a more compact form 
$$
\mathbf{H}_f\, \underline{q} = \underline{b}\,,
$$
where $\underline{q}$ is a vector with the unknown polynomial coefficients.
The entries of the matrix $\mathbf{H}_f$ are given by inner products of scaled monomials defined on the face $f$. 
Therefore, by applying an appropriate quadrature rule for polygons~\cite{Chin:2021:SBC}, they can be computed explicitly. 
The components of the vector $\underline{b}$ correspond exactly to the degrees of freedom of type \texttt{D1} of the face $f$ we are considering scaled by $|f|$, 
and are thus computable even though they involve the integration of the virtual function $v_h$.

\begin{remark}
Given a face $f$ and a virtual basis function $\phi_h$ associated with a degree of freedom of type \texttt{D1}, 
the projection operator is nonzero \emph{only} if $\phi_h$ is associated with the same face $f$. 
\corr{
In particular, the corresponding vector $\underline b$ is different from zero only in the component corresponding to the polynomial that appears in the definition of the degree of freedom of type \texttt{D1}.
}
Moreover, for all basis functions associated with degrees of freedom of type \texttt{D2}, 
the projection operator vanishes.
\end{remark}

\begin{remark}
If we were in two dimensions, 
the faces correspond to edges and one can simply verify that the projection $\Pf$ on an edge does not depend on the edge length. 
\corr{In fact, if we look at the entries of the matrix $\mathbf{H}_f$ and at the definition of the monomials $m^{\jj}$ in one dimension, we can observe that, after a change of variables, the integrals do not depend on the edge size.}
This imply that it is possible to compute the projection one time for all the edges.
Similarly, in three dimensions, it is possible to verify that if two faces have the same shape, the projection does not depend on the face size. 
This is particularly useful in the case where the polyhedrons are formed by the same faces like a cube decomposition of a mesh.
\end{remark}

\subsubsection{The projection $\PN$}\label{sec:projNabla}

Given a polyhedron $E \in \Omega_h$, 
we define a projection operator from the VE space to the polynomial space of order $k$,
$\PN : V_h^k(E) \to \mathbb{P}_k(E)$ based on the $H^1(E)$ semi-inner product.
The operator $\PN$ is such that, for any $p_k \in \mathbb{P}_{k}(E)$, the following condition holds:
\begin{equation}
\bigl(\nabla p_k,\, \nabla(v_h - \PN v_h)\bigr)_E = 0\quad \forall\, p_k \in \mathbb{P}_k(E).
\label{eqn:defPNE}
\end{equation}
As before, $\PN v_h$ is an element of $\mathbb{P}_k(E)$ so  
we can represent such projection operator in the basis $\mathcal{M}^k(E)$:
\begin{equation}
\PN v_h = \sum_{|\jj|\leq k} r^{\jj} m^{\jj}\,.
\label{eqn:pnPolyDef}
\end{equation}
Since Equation~\eqref{eqn:defPNE} involves only the gradients of the polynomials, 
it does not provide enough conditions to uniquely determine all the coefficients of $\PN v_h$
in the chosen polynomial basis $\mathcal{M}^k(E)$.
In particular, the component in the kernel of the gradient operator, 
i.e., the constant part, remains undetermined.
To ensure uniqueness,
we complement the definition by enforcing the following condition, 
which effectively fixes the average value of the projection on the boundary:
\begin{equation}
\int_{\partial E} \PN v_h\,\mathrm{d}f =  \int_{\partial E} v_h\,\mathrm{d}f\,.
\label{eqn:pnFix}
\end{equation}
This choice guarantees the projection is well-posed and computable using only available degrees of freedom.
With a slight abuse of notation, we refer to the integrals in Equation~\eqref{eqn:pnFix} as $(\cdot,1)_{\partial E}$.

Therefore, to find the coefficients $r^{\jj}$, we set the following linear system
\begin{subequations}
\begin{align}
\sum_{|\jj|\leq k} r^{\jj} ( m^{\jj},\,1)_{\partial E} &= (v_h,\,1)_{\partial E} \, , \\
\sum_{|\jj|\leq k} r^{\jj} (\nabla m^{\jj},\,\nabla m^{\ii})_E &= (\nabla v_h,\,\nabla m^{\ii})_E\qquad \forall\, 0 < |\ii|\leq k\,.
\label{eqn:systPN}
\end{align}
\end{subequations}
This linear system can be written in a more compact form as 
\begin{equation}
\mathbf{G}\, \underline{r} = \underline{c}\,.    
\label{eqn:lsnabla}
\end{equation}
The entries of the matrix $\mathbf{G}$ are given by integrals of monomials 
or inner product between gradients of monomials defined inside the polyhedron $E$ so 
they are computable~\cite{Chin:2021:SBC}. 
In order to understand if we are able to compute the unknown coefficients of the vector $\underline{r}$,
we have to check if the right hand side is computable although it contains the virtual function $v_h$.
To compute the first component of the vector $\underline{c}$, 
we observe that
$$
 (v_h,\,1)_{\partial E} = \sum_{f\in\partial E} \int_f v_h\df = \sum_{f\in\partial E} \int_f \Pf v_h\df\,.
$$
\corr{
\begin{remark}
The integral in \eqref{eqn:pnFix} is nonzero only if $v_h$ corresponds to a basis function associated with a degree of freedom of type $\texttt{D1}$, where $m_{\ell} = m_0 = 1$.  
In this case, the integral is equal to the face size.
\end{remark}
}
Then, for all the other components, we integrate by parts and we get
\begin{align*}
(\nabla v_h,\,\nabla m^{\ii}) &= \int_E \nabla v_h \cdot \nabla m^{\ii}\dE\\
&=- \int_E v_h \, \Delta m^{\ii}\dE + \sum_{f \in \partial E} \int_f \left( \nabla m^{\ii} \cdot \nn^E\right)\,v_h\df\,\\
&=- \int_E v_h \, \Delta m^{\ii}\dE + \sum_{f \in \partial E} \int_f \left( \nabla m^{\ii} \cdot \nn^E\right)\,\Pf v_h\df\,.
\end{align*}
Since $\Delta m^{\ii}$ is a polynomial of degree $k-2$,
it is computable from the degrees of freedom $\texttt{D2}$.
Then, since $\nabla m^{\ii} \cdot \nn^E$ is also a polynomial of degree $k-1$,
we can exploit the projection operator introduced in Section~\ref{sec:pfProj},
to compute it.

\begin{remark}
The previous computations were carried out for a generic virtual function $v_h \in V_h^k(E)$.
However, if $v_h$ is taken to be one of the basis functions associated with the degrees of freedom \texttt{D2},
then the boundary integrals vanish identically. 
Conversely, if $v_h$ corresponds to a basis function associated with a boundary degree of freedom \texttt{D1} and 
face $\overline{f}$,
then the volume integral vanishes, and only the boundary contribution on $\overline{f}$ remains.
\end{remark}

\subsubsection{The projection $\PZ$}

Given a polyhedron $E \in \Omega_h$, 
we define a projection operator from the VE space to the polynomial space of order $k$,
$\PZ : V_h^k(E) \to \mathbb{P}_k(E)$ via the $L^2$ inner product:
\begin{equation}
(p_k,\, v_h - \PZ v_h)_E = 0\quad \forall\, p_k \in \mathbb{P}_k(E).
\label{eqn:defPZE}
\end{equation}
As before, the conditions on $ p_k \in \mathbb{P}_k(E)$ can be equivalently replaced by the basis elements $\mathcal{M}^k(E)$
and we can also expand the polynomial $\PZ v_h$ considering the same basis, i.e.,
\begin{equation}
\PZ v_h = \sum_{|\jj|\leq k} s^{\jj} m^{\jj}\,.
\label{eqn:pzPolyDef}
\end{equation}
Hence, to find the coefficients $s^{\jj}$, 
we build the following linear system
\begin{equation}
\sum_{|\jj|\leq k} s^{\jj} (m^{\jj},\, m^{\ii})_E = (v_h,\, m^{\ii})_E \qquad \forall\, 0 \leq |\ii|\leq k\,,
\label{eqn:monPZE}
\end{equation}
which can be written in a more compact form as
\begin{equation}
\mathbf{H}\, \underline{s} = \underline{d}\,.    
\label{eqn:l2Mat}
\end{equation}
The entries of the matrix $\mathbf{H}$ are given by inner product between monomials defined inside the polyhedron $E$,
hence they are computable~\cite{Chin:2021:SBC}. 
The right hand side, $\underline{d}$,  contains the virtual function $v_h$,
and we need to check if these components can be computed
First, we observe that for $0 \leq |\ii|\leq k-2$
\begin{equation}
(v_h,\, m^{\ii})_E = \int_{E} v_h \, m^{\ii}\dE\,,
\label{eqn:cond1PZ}
\end{equation}
we have exactly the degrees of freedom $\texttt{D2}$ so these integrals can be computed.
Finally, when \corr{$k-1 \leq |\ii|  \leq k$},
we can exploit the enhancing property of the space $V^k_h(E)$ to compute them,
\begin{equation}
(v_h,\, m^{\ii})_E = \int_{E} v_h \, m^{\ii}\dE = \int_{E} \PN v_h \, m^{\ii}\dE\,.
\label{eqn:cond2PZ}
\end{equation}

The left hand side of this last equation is an integral of the product of two polynomials 
and thus can be computed using appropriate quadrature rules defined on polyhedra.

\begin{remark}
It is worth noticing that for $k=1$ the projection $\PZ v_h$ is defined \emph{only} with the conditions of Equation~\eqref{eqn:cond2PZ}.
As a consequence, the projection $\PN v_h$ and $\PZ v_h$ coincide.
\end{remark}

\begin{remark}
The previous computations were carried out for a generic virtual function $v_h \in V_h^k(E)$.
If $v_h$ is taken to be one of the basis functions associated with the degrees of freedom \texttt{D1},
Equation~\eqref{eqn:cond1PZ} is always zero.
On the other hand, if it is a basis functions associated with the degrees of freedom \texttt{D2},
only one integral of Equation~\eqref{eqn:cond1PZ} is different from zero and it is equal to the volume of the polyhedron $E$.
\corr{On the other hand, if we need to compute \eqref{eqn:cond2PZ} for a basis function associated with a degree of freedom, we exploit the fact that  
\[
(\PN v_h, m^{\ii})_E = \sum_{|\jj| < k} s^{\jj} (m^{\jj}, m^{\ii})_E \, ,
\]
where $\{s^{\jj}\}$ are the coefficients of $\PN v_h$ obtained in \eqref{eqn:lsnabla}, and $(m^{\jj}, m^{\ii})_E$ corresponds to the column of $\Hmat$ associated with $m^{\ii}$.  
Hence, all these quantities are already available.}
\label{rem:pnCons}
\end{remark}

\subsubsection{The projection $\Pi_{k-1}^{0,E}$}

Given a polyhedron $E \in \Omega_h$, 
we define a projection operator that acts on the gradients of virtual element functions
$\PV: \nabla V_h^k(E) \to [\mathbb P_{k-1}(E)]^3$.
This projection operator is based on the $L^2$ inner product through the condition
\begin{equation}
\left(\pp_{k-1},\, \nabla  v_h - \PV \nabla v_h\right)_E = 0\quad \forall\, \pp_{k-1} \in \left[\mathbb{P}_{k-1}(E)\right]^3.
\label{eqn:defPZV}
\end{equation}
Starting from the scalar monomial basis $\mathcal{M}^{k-1}(E)$,
we construct the following basis for vector evaluated polynomials
$$
\boldsymbol{\mathcal{M}}^{k-1}(E) 
:= \left\{
\begin{bmatrix}
m^{\jj}\\
0\\
0\\
\end{bmatrix},
\begin{bmatrix}
0\\
m^{\jj}\\
0\\
\end{bmatrix},
\begin{bmatrix}
0\\
0\\
m^{\jj}\\
\end{bmatrix}\::\:
0\leq |\jj|\leq k-1
\right\}\,.
$$
The conditions on $\pp_{k-1} \in \left[\mathbb{P}_{k-1}(E)\right]^3$
can be equivalently replaced by the elements of $\boldsymbol{\mathcal{M}}^{k-1}(E)$. 
Hence, the resulting linear system can be written in the form 
\begin{equation}
\sum_{\mm^{\ii}\in\vmb} t^{\jj} (\mm^{\jj},\, \mm^{\ii})_E = (\nabla v_h,\, \mm^{\ii})_E\qquad \forall\, \mm^{\ii}\in\vmb\,,
\label{eqn:monPZGE}
\end{equation}
which can be written in a more compact form as
\begin{equation}
\mathbf{H}_v\, \underline{t} = \underline{e}\,.    
\label{eqn:monPZGEMat}
\end{equation}
\begin{impDet}
Exploiting the structure of the basis elements of $\vmb$,
and ordering its elements so that the first group corresponds to the first component,
the second to the second, and the third to the third one,
the matrix $\mathbf{H}_v$ acquires the following block-diagonal structure
$$
\mathbf{H}_v \coloneqq
\begin{bmatrix}
    \Hmat & 0 & 0 \\
    0 & \Hmat & 0 \\
    0 & 0 & \Hmat
\end{bmatrix} \, ,
$$
where $\Hmat$ is exactly the matrix in Equation~\eqref{eqn:l2Mat} \corr{resized to the polynomials of degree up to $k-1$}.
\end{impDet}
The right hand side contains virtual functions,
hence we have to check if we are able to compute all the integrals.
As for the projection operator $\PN v_h$,
we use the integration by parts and the projection operator $\Pf$:
\begin{eqnarray*}
\int_E \nabla v_h\cdot \mm^{\ii} \dE &=& - \int_E v_h\diver(\mm^{\ii})\dE + \int_{\partial E} (\nn\cdot\mm^{\ii})v_h\df\\
&=& - \int_E v_h\diver(\mm^{\ii})\dE + \sum_{f\in\partial E}\int_f (\nn\cdot\mm^{\ii})\:v_h\df\\
&=& - \int_E v_h\diver(\mm^{\ii})\dE + \sum_{f\in\partial E}\int_f (\nn\cdot\mm^{\ii})\:\Pf v_h\df\,.
\end{eqnarray*}
Since we are considering a vector evaluated polynomial space of order $k-1$,
$\diver{(\mm^{\ii})}$ is a polynomial of degree $k-2$,
and it is computable from the degrees of freedom $\texttt{D2}$.

\begin{remark}
The previous computations were made considering a generic VE function $v_h$
but, if we take a basis VE function,
\corr{these computations simplify in a similar way to Remark~\ref{rem:pnCons}}.
\end{remark}

\begin{impDet}
Since we are considering straight faces, the normal vector is constant,
hence in the computation of the boundary integral only the component of the normal is varying.
In fact, if we consider the following vector evaluated monomial 
$$
\mm_x =\begin{bmatrix}
m^{\ii}\\
0\\
0\\
\end{bmatrix}\,,
$$
we have 
\begin{eqnarray*}
\int_f (\nn_x m^{\ii})\:\Pf v_h\df = \nn_x \int_f m^{\ii}\:\Pf v_h\df\,,
\end{eqnarray*}
where $\nn_x$ is the first component of the normal $\nn$.
For a practical point of view,
it is possible to exploit this fact and save computations.
\end{impDet}

\subsection{Linear and Bilinear Forms}

In this section, we focus on the computation of all the linear and bilinear forms introduced in Section~\ref{sec:model}.
Specifically, since each form admits a decomposition into element-wise contributions, 
we focus on the definition and implementation of the local forms on a generic element $E$.

\subsubsection{The Bilinear form $a(\cdot,\cdot)$}\label{sec:a:const}

In this subsection we consider the bilinear form defined in Equation~\eqref{eq:a-c}.
As usual in the VE framework, the local contributions of this bilinear form are replaced by 
a discrete approximation of them:
$$
a(u,v) = \sum_{E \in \Omega_h} a^E(u,v) \approx \sum_{E \in \Omega_h} a^E_h(u_h,v_h)\,,
$$
where we defined the bilinear form $a_h^E:V_h^k(E) \times V_h^k(E) \to \R$ as 
\begin{multline}
a_h^E(u_h,v_h) 
\coloneqq 
\int_E  \boldsymbol{\kappa} \, \PV \nabla u_h \cdot \ \PV \nabla v_h \, \dE + \\
+\bar{\boldsymbol{\kappa}}_E\StabE_a\bigl( (I-\PN)u_h, (I-\PN)v_h \bigr) \, ,
\label{eq:ahE}
\end{multline}
where $\bar{\boldsymbol{\kappa}}_E$ is the trace of the tensor $\kappa$ 
evaluated at the barycentre of $E$.
The form $\StabE$ is any positive definite bilinear form that satisfies the following properties:
\begin{itemize}
\item it is zero if at least one of the entry is a polynomial of degree $k$,
\item it scales as the continuous form $a^E(\cdot,\cdot)$, i.e., 
there exists two positive values $\alpha_*$ and $\alpha^*$ such that 
\begin{equation}\label{eq:stabE}
\alpha_* a^E( v_h, v_h )
\leq 
\StabE_a( v_h, v_h )
\leq
\alpha^* a^E( v_h, v_h )
\quad
\forall v_h \in \ker(\PN) \,.
\end{equation}
\end{itemize}
There are different choices for the bilinear form $\StabE$,
we refer the reader to~\cite{Mascotto:2023:TRS} for a complete survey.

In Section~\ref{sec:num:kconv}, we make a numerical analysis on the robustness of the method 
with the most common ones: 
the classical \texttt{Dofi-Dofi}~\cite{volley} and 
the \texttt{D-Recipe}~\cite{BDR:2017}.

\begin{remark}
A more classical approach to discretize the bilinear form $ a^E(\cdot,\cdot) $
is considering this approximation 
\begin{multline}\label{eq:ahE-bad}
\tilde a_h^E(u_h,v_h) 
\coloneqq 
\int_E  \boldsymbol{\kappa} \, \nabla \PN u_h \cdot  \nabla \PN v_h \, \mathrm{d}E \\
+
\bar{\boldsymbol{\kappa}}_E\StabE_a\bigl( (I-\PN)u_h,\ (I-\PN)v_h \bigr) \, .
\end{multline}
However, it was shown in \cite{BBMR:2016} for conforming VEM in 2D 
that such a choice results in a significant loss of convergence when $k \geq 3$.
In Section~\ref{sec:num:diffNab}, we give the numerical evidence about this fact for the VE spaces here proposed.
The VE approximation~\eqref{eq:ahE} and~\eqref{eq:ahE-bad} are indeed equivalent for any degree $k$ \emph{only}
in the case where $\boldsymbol{\kappa} \nabla u_h$ is itself a gradient of a function. 
\label{rem:diffGradProj}
\end{remark}

\begin{impDet}
During the assembling the stiffness matrix, 
the integral in Equation~\eqref{eq:ahE} 
can be computed for each polyhedron basis function using a simple matrix-by-matrix product.
Let $\boldsymbol{\Pi}^{0,\text{G}}_*$ be a matrix defined as 
$$
\boldsymbol{\Pi}^{0,\text{G}}_* = \begin{bmatrix}
& & & \\
\underline{t}_1 & \underline{t}_2 & \ldots & \underline{t}_n\\
& & & \\
\end{bmatrix}
$$
whose columns $t_i$ are the solution of the linear system~\eqref{eqn:monPZGEMat} associated with the $i-$th VE basis function, 
and let $\mathbf{H}_{\boldsymbol{\kappa}}$ be a matrix whose entries are 
$$
(\boldsymbol{\kappa} \, \mm^{\jj},\, \mm^{\ii})_E\qquad \forall\, \mm^{\ii},\,\mm^{\jj}\in\vmb\,.
$$
Hence, we can obtain the local matrix associated with the element $E$ 
that represent the first integral of Equation~\eqref{eq:ahE} 
via this chain of matrix product
$$
\left(\boldsymbol{\Pi}^{0,\text{G}}_*\right)^{\top}\,\mathbf{H}_{\boldsymbol{\kappa}}\:\boldsymbol{\Pi}^{0,\text{G}}_*\,.
$$
\label{impDet:grad}
\end{impDet} 

\subsubsection{The Bilinear form $m(\cdot,\cdot)$}

The procedure for constructing the bilinear form $m(\cdot,\cdot)$ is similar to the one described in Section~\ref{sec:a:const}.
Similarly, the continuous form is decomposed over the polyhedra and
replaced by a discrete approximation:
$$
m(u,v) = \sum_{E \in \Omega_h} m^E(u,v) \approx \sum_{E \in \Omega_h} m^E_h(u_h,v_h)\,,
$$
where we defined the bilinear form $m_h^E:V_h^k(E) \times V_h^k(E) \to \R$ as 
\begin{equation}\label{eq:mhE}
m_h^E(u_h,v_h) 
\coloneqq 
\int_E \PZ u_h \, \PZ v_h \, {\rm d }E 
+
\StabE_m \bigl( (I-\PZ)u_h, (I-\PZ)v_h \bigr) \,.
\end{equation}
The stabilization term $\StabE_m(\cdot,\cdot)$ \corr{is a positive definite bi-linear form that
vanishes} if at least one of the entries is a polynomial of degree $k$,
and it must scale consistently with the continuous form $m^E(\cdot,\cdot)$, i.e., 
\begin{equation}\label{eq:stabBE}
\mu_* m^E(v_h,\,v_h)
\leq 
\StabE_m( v_h, v_h )
\leq
\mu^* m^E(v_h,\,v_h) 
\quad
\forall v_h \in \ker(\PZ) \, .
\end{equation}

\begin{impDet}
During the assembly of the mass term of the linear system, 
the integral in Equation~\eqref{eq:mhE} can be efficiently computed for each polyhedron basis function 
via a simple matrix-by-matrix product as its counterpart in Equation~\eqref{eq:ahE}.
Let $\boldsymbol{\Pi}^0_*$ be a matrix defined as 
$$
\boldsymbol{\Pi}^0_* = \begin{bmatrix}
& & & \\
\underline{s}_1 & \underline{s}_2 & \ldots & \underline{s}_n\\
& & & \\
\end{bmatrix}
$$
whose columns $s_i$ are the solutions 
of the linear system~\eqref{eqn:l2Mat} associated with the $i-$th VE basis function.
Hence, we can obtain the local matrix associated with the element $E$, 
representing the first integral of Equation~\eqref{eq:mhE}, 
can be obtained via the following chain of matrix products:
$$
\left(\boldsymbol{\Pi}^0_*\right)^{\top}\,\mathbf{H}\:\boldsymbol{\Pi}^0_*\,,
$$
where $\mathbf{H}$ is the \emph{same} matrix as the linear system~\eqref{eqn:l2Mat}.
\label{impDet:mass}
\end{impDet} 

\begin{impDet}
Considering the classical \texttt{Dofi-Dofi} stabilization,
both stabilization terms $\StabE_a(\cdot,\cdot)$ and $\StabE_m(\cdot,\cdot)$ can be expressed as a sequence of matrix products~\cite{autostoppisti}.
Let $\mathbf{D}$ be the matrix defined as follows
$$
\mathbf{D} = 
\begin{bmatrix}
dof_1(m^{\jj_1}) & dof_1(m^{\jj_2}) &\ldots &dof_1(m^{\jj_{n_c}})\\
dof_2(m^{\jj_1}) & dof_2(m^{\jj_2}) &\ldots &dof_2(m^{\jj_{n_c}})\\
\vdots           &\vdots            &\ddots &\vdots \\
dof_{N_E}(m^{\jj_1}) & dof_{N_E}(m^{\jj_2}) &\ldots &dof_{N_E}(m^{\jj_{n_c}})
\end{bmatrix}\,.
$$
where $dof_i(\cdot)$ is a linear functional representing the $i-$th degree of freedom.
The stabilization $\StabE_a(\cdot,\cdot)$ becomes
\begin{equation}\label{eq:dofia_mat}
h_E\left(\mathbf{I}-\mathbf{D}\,\boldsymbol{\Pi}^{\nabla}_*\right)^{\top}\,
\left(\mathbf{I}-\mathbf{D}\,\boldsymbol{\Pi}^{\nabla}_*\right)\,,
\end{equation}
where $\boldsymbol{\Pi}^{\nabla}_*$ is a matrix whose columns are the solution of the linear system~\eqref{eqn:lsnabla},
and $\mathbf{I}$ is the identity matrix.
Similarly, the stabilization $\StabE_m(\cdot,\cdot)$ is given by
\begin{equation}\label{eq:dofim_mat}
|E|\left(\mathbf{I}-\mathbf{D}\,\boldsymbol{\Pi}^0_*\right)^{\top}\,
\left(\mathbf{I}-\mathbf{D}\,\boldsymbol{\Pi}^0_*\right)\,,
\end{equation}
where $\boldsymbol{\Pi}^0_*$ is defined in Implementation Detail~\ref{impDet:mass}.

\corr{
The so-called \texttt{D-Recipe} stabilization, introduced in~\cite{BDR:2017}, 
was proposed as an alternative to the \texttt{Dofi-Dofi} term, 
which may become too restrictive for large~$k$, particularly in three-dimensional settings.
In this approach, different basis functions may lead to different energy values. 
To account for this, the following diagonal stabilization form is adopted:
\[
\tilde{\StabE}_a(\varphi_i,\varphi_i) \coloneqq 
\max \{ h_E, (\PV \nabla \varphi_i, \PV \nabla \varphi_i)_E \},
\]
where the set~$\{\varphi_i\}$ denotes the nodal basis functions of~$V_h^k(E)$.
Accordingly, the matrix \eqref{eq:dofia_mat} (and analogously \eqref{eq:dofim_mat}) should be replaced by
\[
\left(\mathbf{I}-\mathbf{D}\,\boldsymbol{\Pi}^{\nabla}_*\right)^{\top}\,
\tilde{\mathbf{R}} \, 
\left(\mathbf{I}-\mathbf{D}\,\boldsymbol{\Pi}^{\nabla}_*\right)\,,
\]
where $\tilde{\mathbf{R}}$ is a diagonal matrix defined by
\[
(\tilde{\mathbf{R}})_{ii} 
= 
\max \{ h_E, (\PV \nabla \varphi_i, \PV \nabla \varphi_i)_E \}
\, .
\]
We conclude by observing that the quantity $(\PV \nabla \varphi_i, \PV \nabla \varphi_i)_E$ is already required in the assembly of the local stiffness matrix. Therefore, no additional computations are needed to construct the stabilization matrix.
}%End corr
\end{impDet}

\subsubsection{The Right-Hand Side}

\corr{The load term is also split among the contributions of each polyhedron}, 
and the virtual function is simply replaced with its $L^2$ projection:
$$
F(v)=\sum_{E\in\Omega_h} \int_E f\,v\dE\approx \sum_{E\in\Omega_h} F_h^E(v_h)\,,
$$
where we define the linear form $F_h^E:V_h^k(E) \to \R$ as
\begin{equation}
F_h^E(v_h) := \int_E f\,\PZ v_h\dE\,.
\label{eqn:forceE}
\end{equation}
It is worth noting that in this paper we consider a more accurate approximation of the load term 
than those proposed in~\cite{blanca,coreani}.
This choice does not improve the convergence rate but
marginally reduces the error of the discrete solution.
As shown in~\cite{blanca,coreani},
using a lower order approximation does not affect the convergence rate of the error measured in the $L^2$ norm or in the $H^1$ seminorm.

\begin{impDet}
When assembling the Right-Hand Side of the linear system, 
the integral in Equation~\eqref{eqn:forceE} can be computed for each polyhedron basis function using a simple matrix-vector product.
Let $\mathbf{H}_{f}$ be a vector whose entries are 
$$
(f,\, m^{\ii})_E\qquad \forall\, m^{\ii},\forall\, 0 \leq |\ii|\leq k\,.
$$
Hence, we can obtain the local vector associated with the element $E$ via the following matrix-vector product:
$$
(\boldsymbol{\Pi}^0_*)^\top \mathbf{H}_{f}\,,
$$
where $\boldsymbol{\Pi}^0_*$ is the matrix introduced in Implementation Detail~\ref{impDet:mass}.
\end{impDet}

\subsection{Discrete Problem}

Based on the definitions of the global spaces, projection operators, 
and the linear and bilinear forms, we can now introduce the discrete problem:
\begin{equation}
\left\{
\begin{aligned}
&\text{find $u_h \in V_h^k(\Omega_h)$ such that:} \\
&a_h(u_h,v_h) + m_h(u_h,v_h) = F_h(v_h) \quad \forall v_h \in V_h^k(\Omega_h) \, ,
\end{aligned}
\right.
\label{eq:discreto}
\end{equation}
where the global bilinear forms are defined as
$$
a_h(u_h,v_h) \coloneqq \sum_{E\in\Omega_h}a_h^E(u_h,v_h)\qquad\text{and}\qquad
m_h(u_h,v_h) \coloneqq \sum_{E\in\Omega_h}m_h^E(u_h,v_h)\,,
$$
and the linear form 
$$
F_h(v_h) \coloneqq \sum_{E\in \Omega_h} F_h^E(v_h)\,.
$$

%% file: dofFig.tex
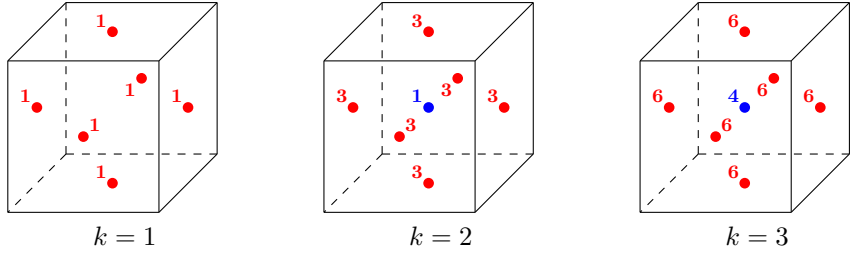
\begin{figure}
\begin{minipage}{0.31\textwidth}
\begin{center}
\begin{tikzpicture}
    \drawcube
    \redcircles{1}
\end{tikzpicture}
\newline
$k=1$
\end{center}
\end{minipage}
\hfill
\begin{minipage}{0.31\textwidth}
\begin{center}
\begin{tikzpicture}
    \drawcube
    \redcircles{3}
    \bluecircle{1}
\end{tikzpicture}
\newline
$k=2$
\end{center}
\end{minipage}
\hfill
\begin{minipage}{0.31\textwidth}
\begin{center}
\begin{tikzpicture}
    \drawcube
    \redcircles{6}
    \bluecircle{4}
\end{tikzpicture}
\newline
$k=3$
\end{center}
\end{minipage}
\caption{Degrees of freedom for a cube. In red we have the Dofs associated to a face of the polyhedron, while in blue we have the moments associated to the interior of the polyhedron.}\label{fig:dofs}
\end{figure}

%% file: numExe.tex
\section{Numerical Experiments}\label{sec:numExe}

In this section, we investigate the numerical behaviour of the nonconforming scheme in different settings.
First, we conduct a convergence analysis with respect to the polynomial degree $k$, 
considering different choices for the stabilization term (see Section~\ref{sec:num:kconv}), 
and a convergence analysis with respect to the mesh size $h$ (see Section~\ref{sec:num:hconv}).
In Section~\ref{sec:num:diffNab}, we provide numerical evidence supporting Remark~\ref{rem:diffGradProj}.
Finally, in Section~\ref{sec:num:conf}, we compare the conforming and nonconforming VEM from an algebraic view point, i.e.,
we analyse execution time, dimension of the linear system, and the sparsity of the global matrix.
All the numerical experiments are implemented using the \texttt{C++} library \texttt{vem++}~\cite{vem++}.

\begin{figure}[!htb]
\begin{center}
\begin{tabular}{ccc}
\includegraphics[width=0.28\textwidth]{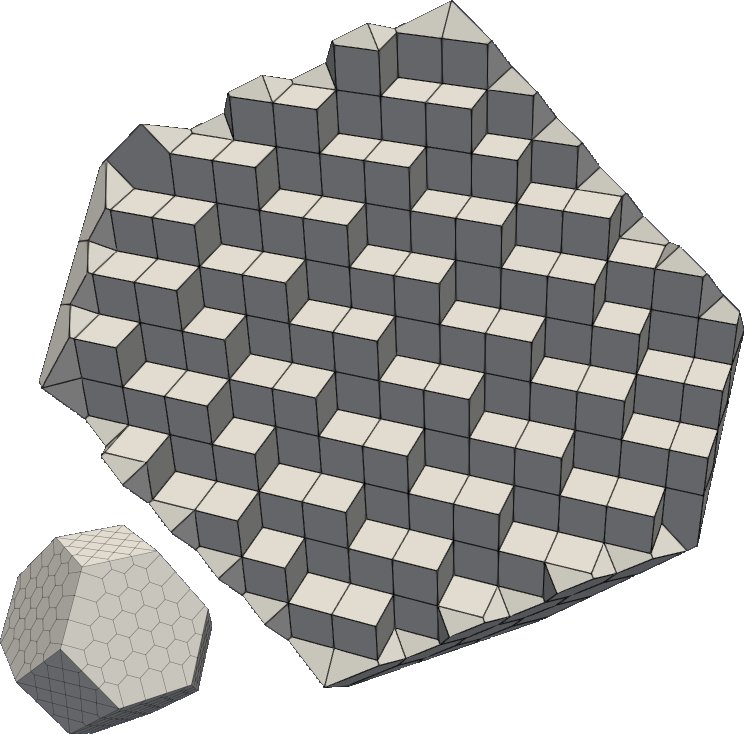} &
\includegraphics[width=0.28\textwidth]{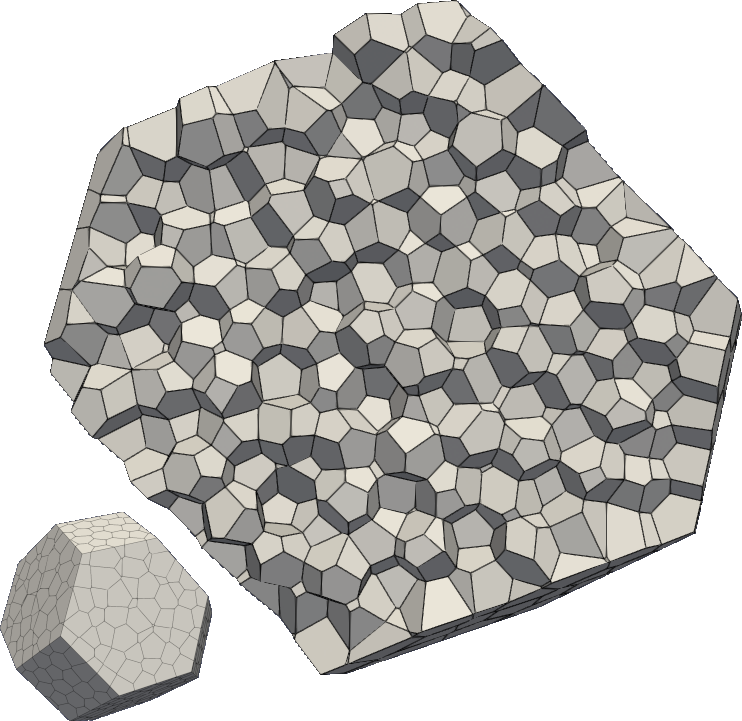} &
\includegraphics[width=0.28\textwidth]{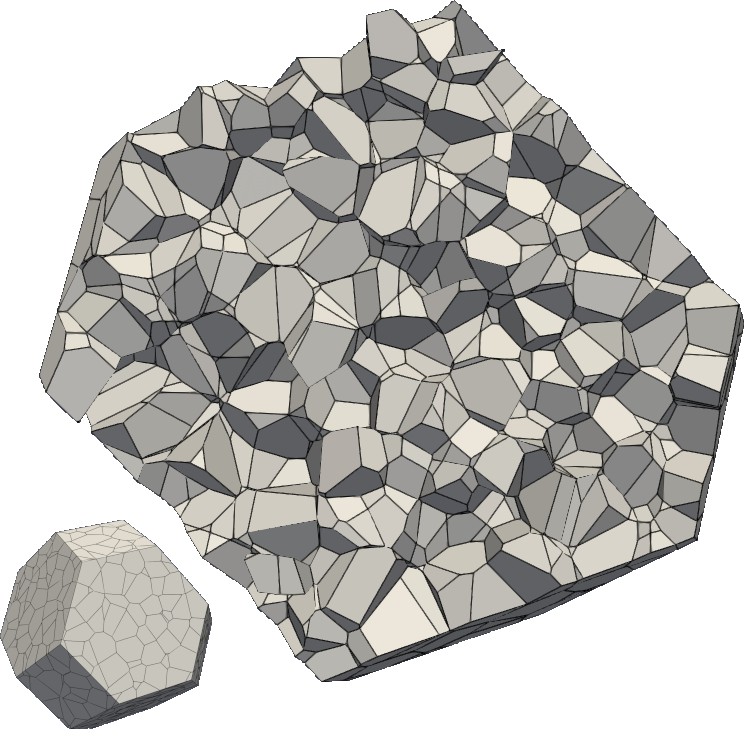} \\
\end{tabular}
\end{center}
\caption{The interior of each type of mesh where the domain $\Omega$ is a truncated octahedron.}
\label{fig:inside}
\end{figure}

In all these experiments,
we consider the following error indicators 
and numerically verify the expected theoretical convergence rate 
\begin{itemize}
\item \textbf{$H^1$-seminorm error}
\[
e_{H^1} := \sqrt{\sum_{E\in\Omega_h} \left\|\nabla(u - \PN u_h)\right\|^2_{0,E}}\approx h^k,
\]
\item \textbf{$L^2$-norm error}
\[
e_{L^2} := \sqrt{\sum_{E\in\Omega_h} \left\|u - \PZ u_h\right\|^2_{0,E}}\approx h^{k+1}.
\]
\end{itemize}
To this end,
we consider the problem defined in Equation~\ref{eqn:strong},
where the domain $\Omega$ is the truncated octahedron,
$\boldsymbol{\kappa}$ is the identity tensor and  
the right-hand side $f$ is chosen so that the exact solution is  
\begin{equation}
u(x, y, z) \coloneqq \sin(\pi x) \cos(\pi y) \cos(\pi z)\,,
\label{eqn:solRefH}
\end{equation}
and we construct a sequence of four meshes with decreasing mesh size $h$.
In order to evaluate the robustness of the nonconforming VEM with respect to element distortion,
we consider three different mesh types of Voronoi tessellations:
\begin{itemize}
    \item \texttt{stru}: seed points are uniformly distributed inside the domain,
    \item \texttt{voro}: the positions of the seed points is optimized via a Lloyd algorithm~\cite{Du:1999:CVT},
    \item \texttt{rand}: seed points are randomly distributed inside the domain.
\end{itemize}
These three types of tessellation are characterised by an increasing level of distortion of the polyhedra.
Indeed, \texttt{stru} meshes have almost hexahedral elements, 
\texttt{voro} meshes may contain polyhedra with small faces and/or edges, 
and \texttt{rand} meshes typically exhibit highly distorted elements.
In Figure~\ref{fig:inside},
we show the interior of these mesh types to
have an idea of the distortion of each mesh type.
In the following examples, 
we refer to these meshes using their type and an integer number from 1 to 4.
Hence, \texttt{rand2} refers to the second refinement level of the domain using seed points randomly distributed inside the domain. 

\subsection{Convergence analysis in $k$}\label{sec:num:kconv}

In this subsection, we test the convergence with respect to the approximation degree $k$,  
considering two choices for the stabilization term $\StabE_a$:  
the \texttt{Dofi-Dofi} stabilization introduced in~\cite{volley} and  
the \texttt{D-Recipe} proposed in~\cite{BDR:2017}.  
According to Remark 6.3 in~\cite{autostoppisti},  
the stabilization term $\StabE_m$ is only required when the equation is reaction-dominated,  
and thus can be set to zero in this case.

\begin{figure}[!htb]
\begin{center}
\begin{tabular}{ccc}
\includegraphics[width=0.44\textwidth]{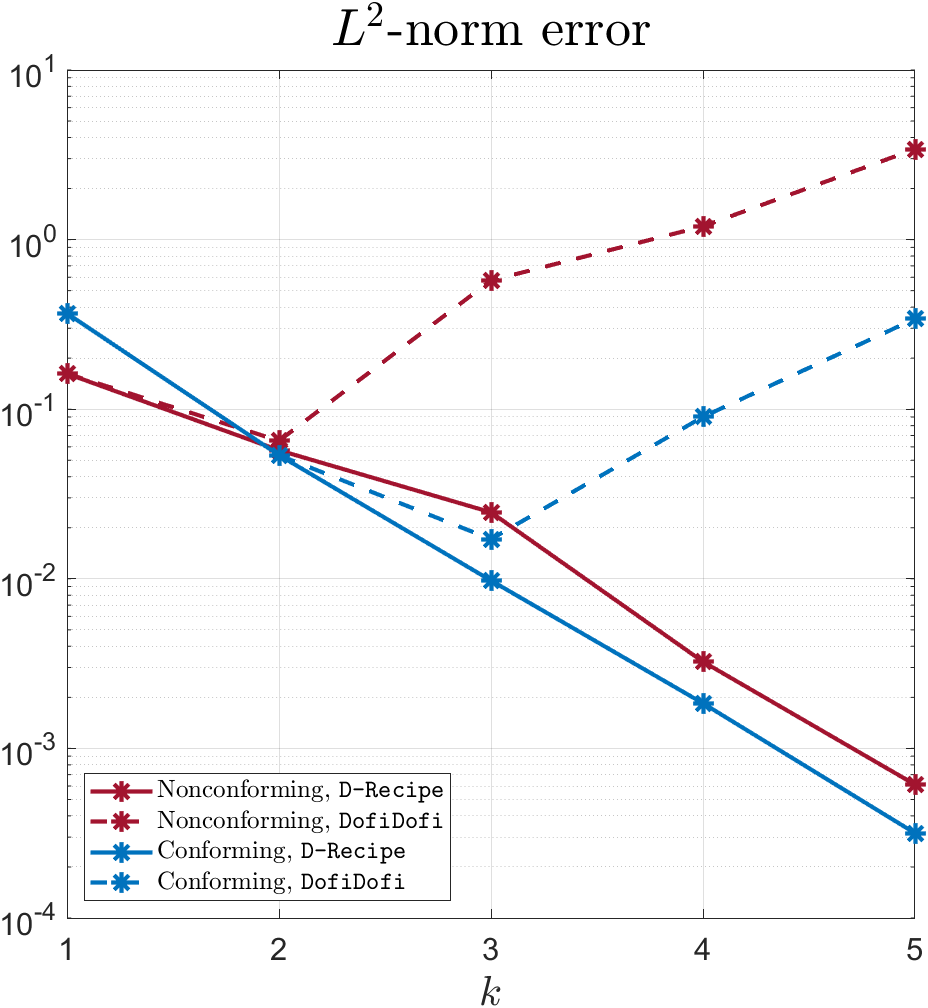} &\hfill&
\includegraphics[width=0.44\textwidth]{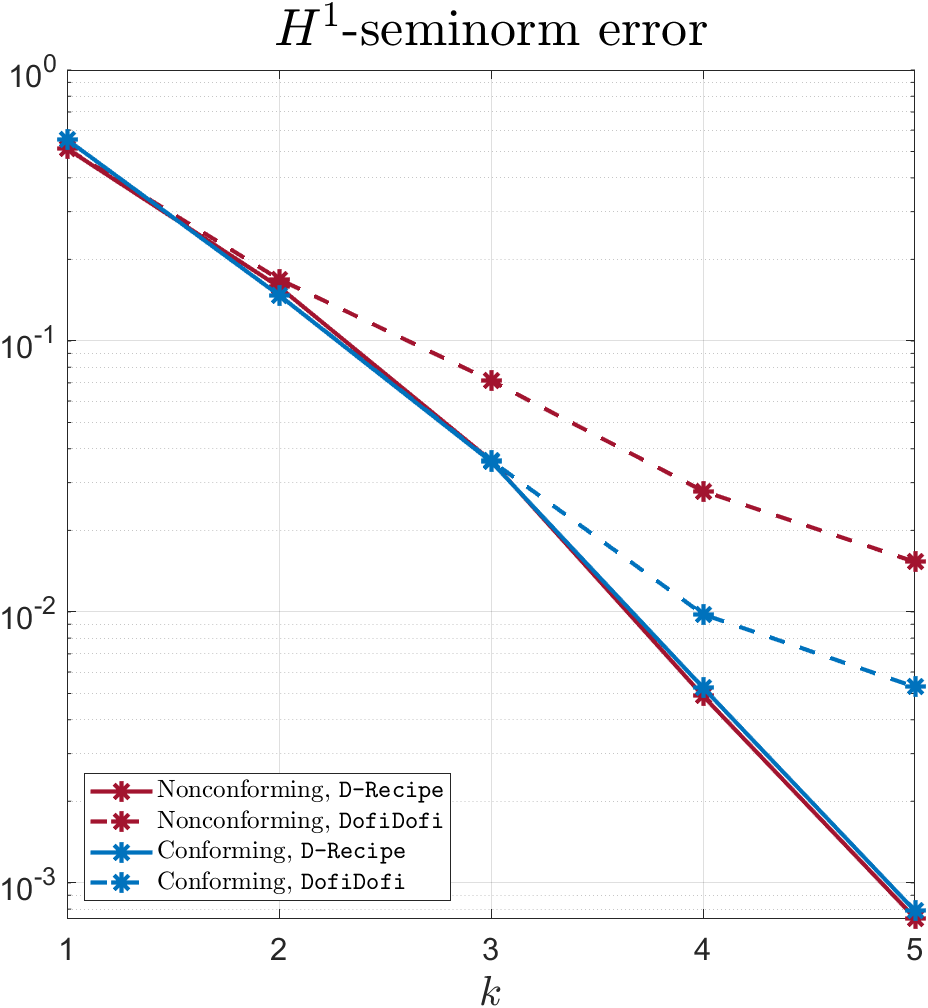} \\
\end{tabular}
\end{center}
\caption{Convergence analysis in $k$: the trends of $e_{H^1}$ and $e_{L^2}$ as the approximation degree $k$ varies.}
\label{fig:exe1:k}
\end{figure}

We fix the mesh \texttt{voro2} and show in Figure~\ref{fig:exe1:k} the convergence rates of both $e_{H^1}$ and $e_{L^2}$.  
Similar results have been observed for the meshes \texttt{stru2} and \texttt{rand2}.

The results show exponential convergence in terms of $k$ when using the \texttt{D-Recipe} stabilization,  
whereas both errors exhibit a slight increase for $k = 4$ and $k = 5$ when using the \texttt{Dofi-Dofi} stabilization.  
These results are consistent with those reported in~\cite{BDR:2017} for the conforming case (see Section 3.4 therein).

Furthermore, to better highlight this similar behaviour,  
we have conducted the same experiments using conforming VE spaces as well; 
see the blue lines in Figure~\ref{fig:exe1:k}.

\subsection{Convergence analysis in $h$}\label{sec:num:hconv}

In this section, we analyse the convergence rate for the nonconforming VEM for $k = 1$, 2, and 3,  
comparing the results obtained with the mesh types \texttt{stru}, \texttt{voro}, and \texttt{rand}.  
We employ the \texttt{D-Recipe} stabilization term,  
since the previous example has shown  
that this choice provides a more accurate trend for higher degrees.

\begin{figure}[!htb]
\begin{center}
\begin{tabular}{ccc}
\includegraphics[width=0.44\textwidth]{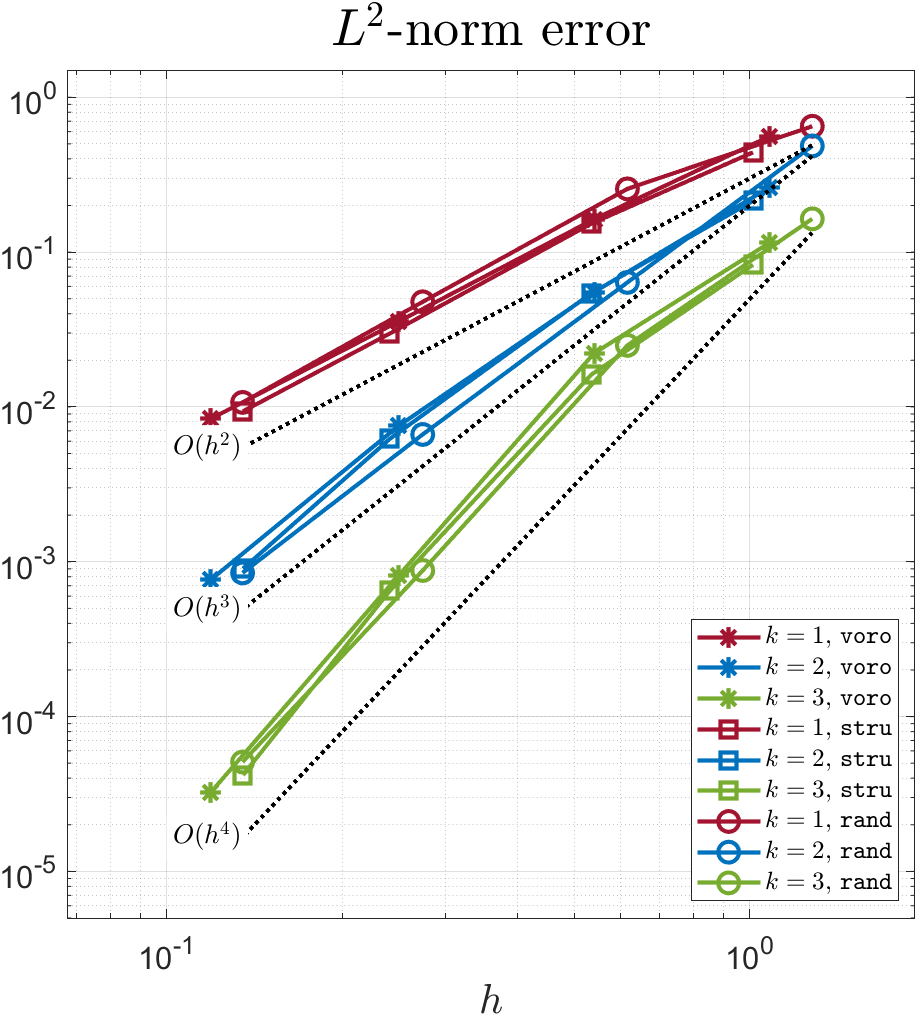} &\hfill&
\includegraphics[width=0.44\textwidth]{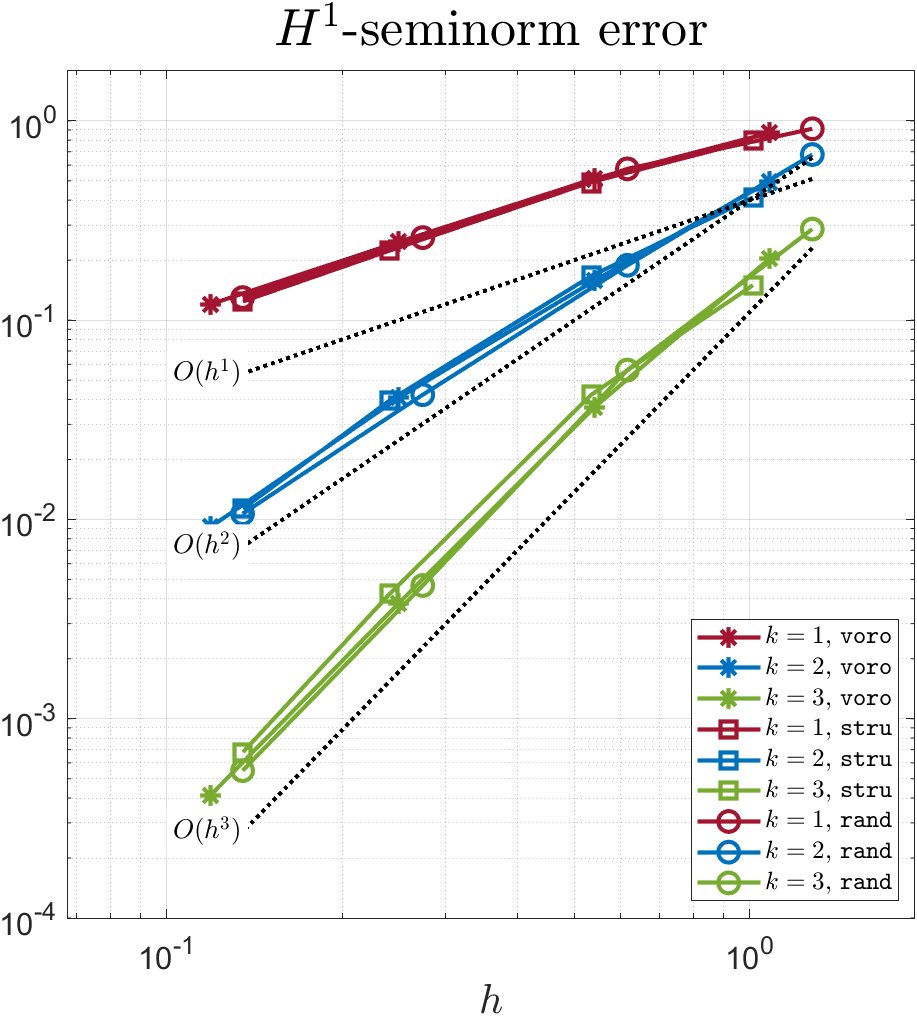} \\
\end{tabular}
\end{center}
\caption{Convergence analysis in $h$:  the trends of $e_{H^1}$ and $e_{L^2}$ varying $h$ and mesh types.}
\label{fig:exe1}
\end{figure}

In Figure~\ref{fig:exe1}, we show the error trends.  
The convergence lines follow the expected theoretical decay:  
$e_{L^2}$ and $e_{H^1}$ converge as $h^{k+1}$ and $h^k$, respectively.  

Furthermore, the convergence lines associated with the different mesh types are close to each other,  
indicating that the proposed method is robust with respect to element distortion.

%
% DA FARE: setta a 0 la stabilizzazione sulla massa, non serve in questo caso 
% mesh  quelle in "../inputMesh/patagon/patagon_random, ../inputMesh/patagon/patagon_regular ../inputMesh/patagon/patagon_structured,"
% k=1,2,3 ,
% nel grafico: - ogni grado ha un colore differente 
%              - le mesh types hanno un marker differente
%
% a Dio piacendo le linee associate ai gradi dovrebbero essere vicine quindi lascerei il commento rosso
%

\subsection{Effects of the discretization of the diffusion term}\label{sec:num:diffNab}

In this section, we compare the discretizations~\eqref{eq:ahE} and~\eqref{eq:ahE-bad}.
As already stated in Remark~\ref{rem:diffGradProj}, 
if the bilinear form $a(\cdot,\cdot)$ involves a non-constant diffusion tensor,
the approximation defined in Equation~\eqref{eq:ahE-bad} does not achieve the optimal convergence rate,
whereas the discretization in~\eqref{eq:ahE} behaves as expected.

To this end, we modify the right-hand side $f$  
so that the exact solution remains the one defined in Equation~\eqref{eqn:solRefH}, but now
\[
\boldsymbol{\kappa}(x,y,z)
\coloneqq
\begin{bmatrix}
y^2 + 1 & -xy & -xz \\
-xy & z^2 + 1 & -yz \\
-xz & -yz & x^2 + 1
\end{bmatrix} \,.
\]

We analyse the convergence rate of the nonconforming VEM for $k = 1, 2, 3$, and 4,  
using the \texttt{D-Recipe} stabilization and considering only the \texttt{voro} mesh type.
Similar results were observed for the other mesh types.

\begin{figure}[!htb]
\begin{center}
\begin{tabular}{ccc}
\includegraphics[width=0.44\textwidth]{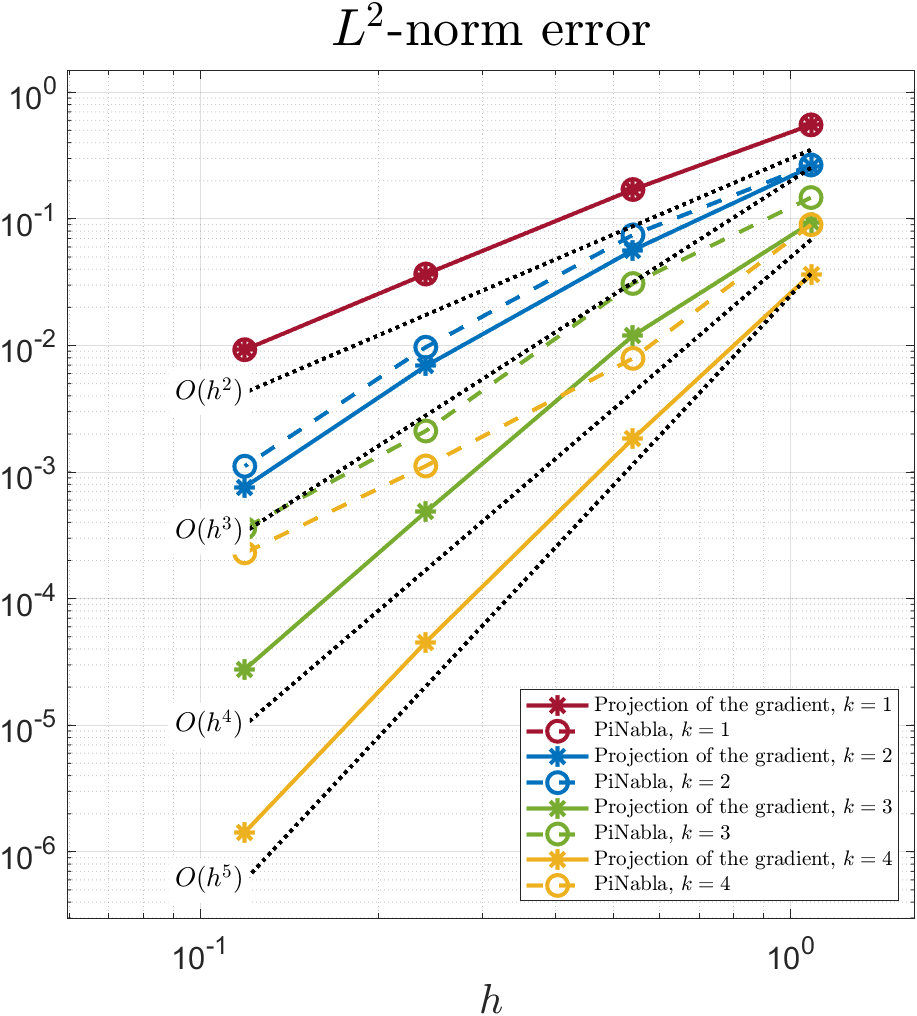} &\hfill&
\includegraphics[width=0.44\textwidth]{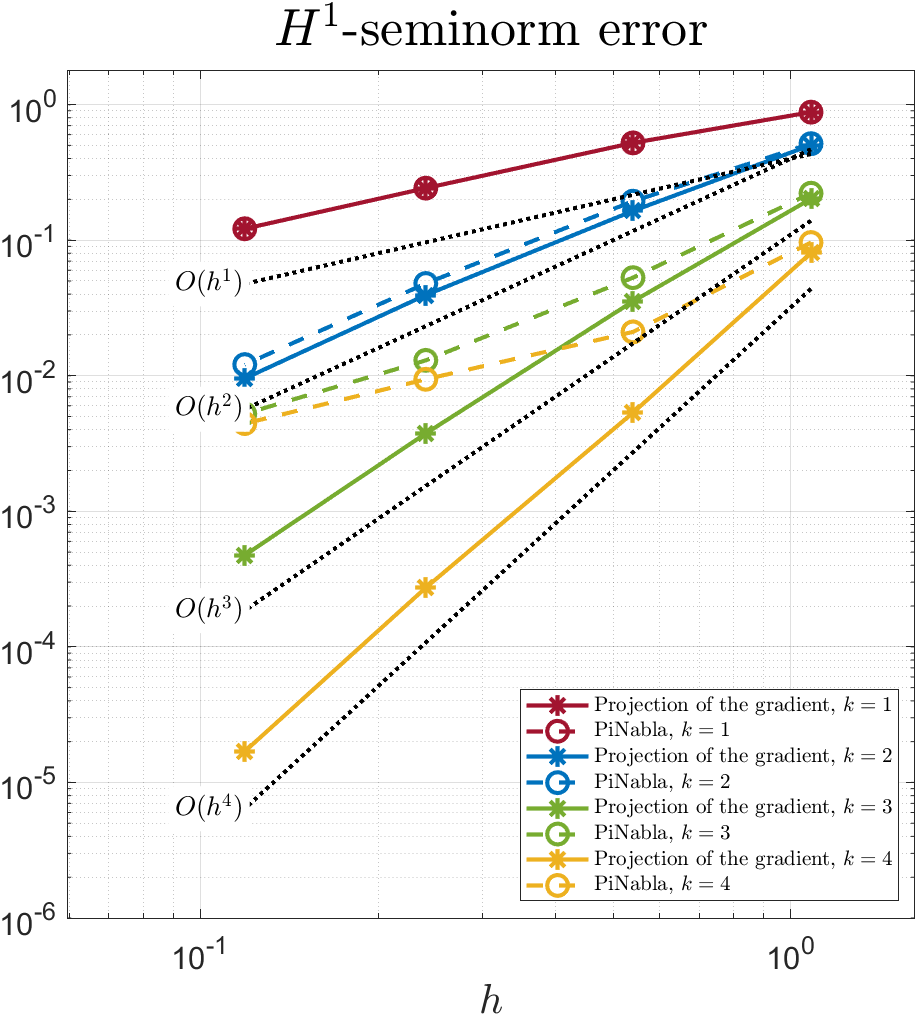} \\
\end{tabular}
\end{center}
\caption{Effects of the discretization of the diffusion term:
comparison between the usage of \eqref{eq:ahE-bad} and \eqref{eq:ahE}.}
\label{fig:exe4}
\end{figure}

In Figure~\ref{fig:exe4}, we show the convergence lines.  
The behaviour of the nonconforming VEM in 3D is consistent with the conforming case in 2D:  
there is a significant loss of convergence rate for $k \geq 3$ when using the discretization~\eqref{eq:ahE-bad}.
However, we also observe that for $k=1$, the two methods produce the same results. 
This makes the choices in \eqref{eq:ahE-bad} preferable in the case $k=1$, 
since they require one fewer projection. 
This result is in line with the analysis in \cite{BBMR:2016}. 
For $k=2$, the two methods behave similarly, 
but the errors are not exactly the same, although they are comparable.

%
% DA FARE: setta a 0 la stabilizzazione sulla massa, non serve in questo caso 
% mesh  quelle in "../inputMesh/patagon/patagon_voro, "
% k=1,2,3 4,
% nel grafico: - ogni grado ha un colore differente 
%              - le mesh types hanno un marker differente
%              - la stab brutta la facciamo tratteggiata quella bella non tratteggiata
%
% direi di fare anche i casi k=1 e 2 per complettezza
%

\subsection{Comparison with the Conforming Method}\label{sec:num:conf}

In this section, we analyse in more detail the characteristic of the global matrix 
arising from the VE discretization of the model problem~\eqref{eqn:strong} and
the execution time of the solver used.

To this end, 
we focus on the meshes \texttt{voro1}, \texttt{voro2} and \texttt{voro3}.
and the direct solver \texttt{MUMPS}~\cite{Amestoy:2001:AFA}
which is integrated in \texttt{vem++} via the \texttt{Petsc} library~\cite{petsc}.
Similar results were observed with other types of meshes.

\begin{table}[h]
\begin{tabular}{|ccccccr|}
\multicolumn{1}{c}{} & \multicolumn{3}{c}{Conforming}                                                                 & \multicolumn{3}{c}{Nonconforming}             \\
\hline
\multicolumn{1}{|c|}{Elements}                          & \multicolumn{1}{c|}{\#Dofs} & \multicolumn{1}{c|}{\%Nonzero}  & \multicolumn{1}{c|}{Time}       & \multicolumn{1}{c|}{\#Dofs} & \multicolumn{1}{c|}{\%Nonzero}  & Time      \\ \hline
\multicolumn{7}{|c|}{$k=1$}                                                                                                                                                                                                     \\ \hline
\multicolumn{1}{|r|}{20}                        & \multicolumn{1}{r|}{41}     & \multicolumn{1}{r|}{62.04\%}  & \multicolumn{1}{r|}{0.102}   & \multicolumn{1}{r|}{77}     & \multicolumn{1}{r|}{20.02\%}  & 0.083  \\ \hline
\multicolumn{1}{|r|}{116}                       & \multicolumn{1}{r|}{404}    & \multicolumn{1}{r|}{12.15\%}  & \multicolumn{1}{r|}{0.454}   & \multicolumn{1}{r|}{597}    & \multicolumn{1}{r|}{3.53\%}  & 0.392  \\ \hline
\multicolumn{1}{|r|}{1001}                      & \multicolumn{1}{r|}{4643}   & \multicolumn{1}{r|}{1.31\%}  & \multicolumn{1}{r|}{4.364}   & \multicolumn{1}{r|}{6064}   & \multicolumn{1}{r|}{0.39\%} & 3.677  \\ \hline
\multicolumn{7}{|c|}{$k=2$}                                                                                                                                                                                                     \\ \hline
\multicolumn{1}{|r|}{20}                        & \multicolumn{1}{r|}{237}    & \multicolumn{1}{r|}{32.94\%}  & \multicolumn{1}{r|}{0.498}   & \multicolumn{1}{r|}{251}    & \multicolumn{1}{r|}{18.45\%}  & 0.344  \\ \hline
\multicolumn{1}{|r|}{116}                       & \multicolumn{1}{r|}{2003}   & \multicolumn{1}{r|}{6.50\%}  & \multicolumn{1}{r|}{2.599}   & \multicolumn{1}{r|}{1907}   & \multicolumn{1}{r|}{3.32\%}  & 1.820  \\ \hline
\multicolumn{1}{|r|}{1001}                      & \multicolumn{1}{r|}{21415}  & \multicolumn{1}{r|}{0.74\%} & \multicolumn{1}{r|}{29.457}  & \multicolumn{1}{r|}{19193}  & \multicolumn{1}{r|}{0.37\%} & 17.830 \\ \hline
\multicolumn{7}{|c|}{$k=3$}                                                                                                                                                                                                     \\ \hline
\multicolumn{1}{|r|}{20}                        & \multicolumn{1}{r|}{550}    & \multicolumn{1}{r|}{26.49\%}  & \multicolumn{1}{r|}{1.964}   & \multicolumn{1}{r|}{542}    & \multicolumn{1}{r|}{17.17\%}  & 1.555  \\ \hline
\multicolumn{1}{|r|}{116}                       & \multicolumn{1}{r|}{4431}   & \multicolumn{1}{r|}{5.23\%}  & \multicolumn{1}{r|}{10.268}  & \multicolumn{1}{r|}{4046}   & \multicolumn{1}{r|}{3.13\%}  & 8.240  \\ \hline
\multicolumn{1}{|r|}{1001}                      & \multicolumn{1}{r|}{46253}  & \multicolumn{1}{r|}{0.60\%} & \multicolumn{1}{r|}{123.262} & \multicolumn{1}{r|}{40388}  & \multicolumn{1}{r|}{0.36\%} & 83.053 \\ \hline
\end{tabular}
\label{tb:exe4}
\caption{Comparison with the conforming method: matrix statistics.}
\end{table}

The data are reported in Table~\ref{tb:exe4}.
The nonconforming method requires less execution time, 
and this advantage becomes more evident increasing the order $k$ and/or decreasing $h$.
It is worth noticing that we take the execution time provided by \texttt{MUMPS},
so the comparison is fair and does not depend on the type of the method but only on the global matrix itself.

\begin{figure}[!htb]
\begin{center}
\begin{tabular}{ccc}
\multicolumn{3}{c}{conforming VE}\\
$k=1$ & &$k=2$\\
\includegraphics[width=0.44\textwidth]{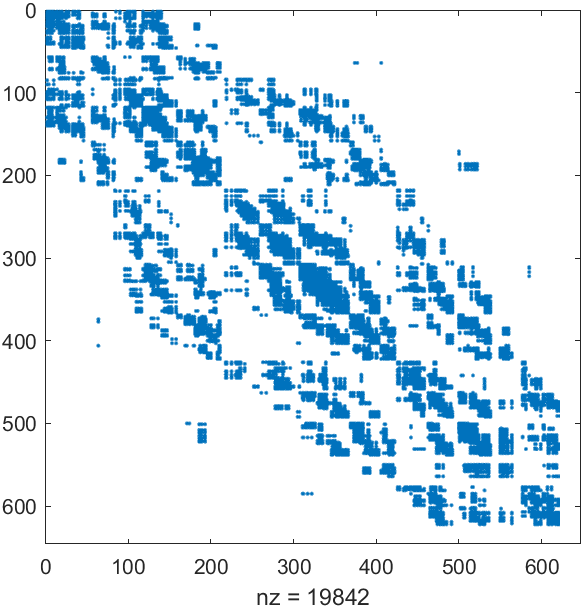} &\hfill&
\includegraphics[width=0.44\textwidth]{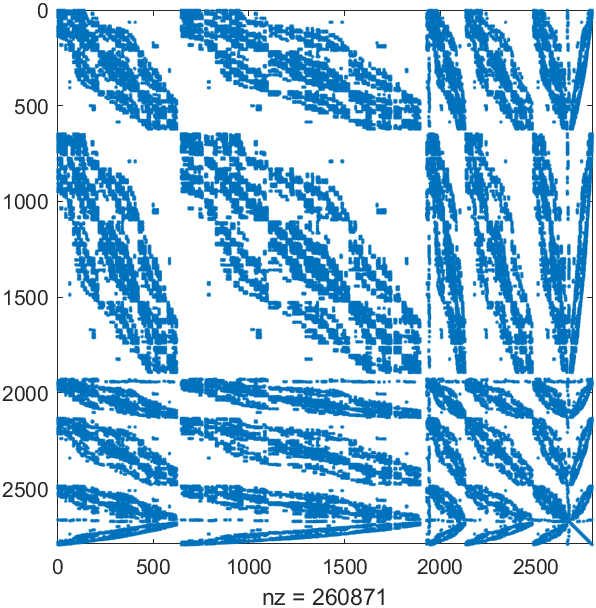} \\
\multicolumn{3}{c}{nonconforming VE}\\
$k=1$ & &$k=2$\\
\includegraphics[width=0.44\textwidth]{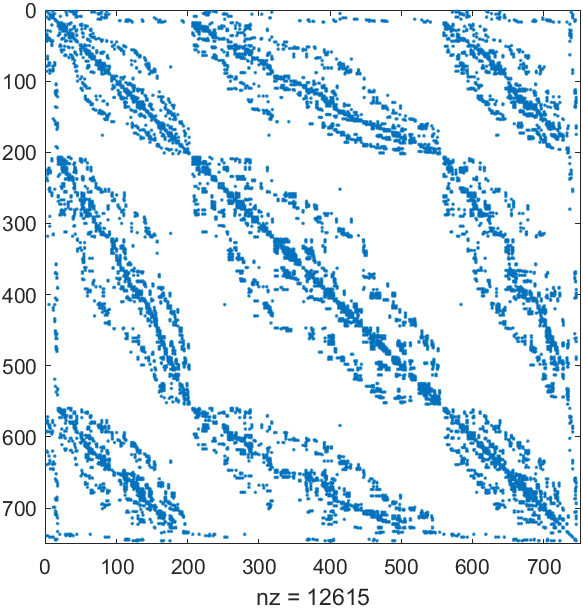} &\hfill&
\includegraphics[width=0.44\textwidth]{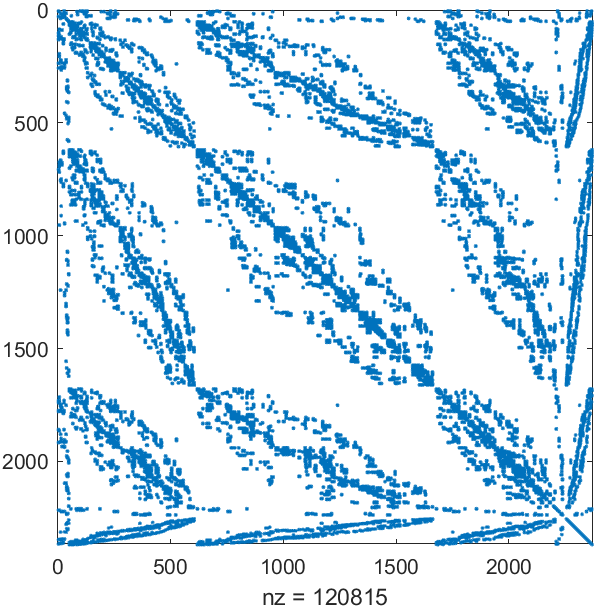} 
\end{tabular}
\end{center}
\caption{Comparison with the conforming method: sparsity pattern of the global matrixes 
obtained from the mesh \texttt{voro2} using VE approximations of degree $k=1$ and $k=2$.}
\label{fig:spy}
\end{figure}

Regarding the structure of the global matrix,
the data shows that the one arising from the nonconforming VEM is sparser than the conforming one.
The difference in sparsity decreases as the polynomial order $k$ increases.
This behaviour is expected, given the locality of interactions in the nonconforming formulation.
Indeed, in the nonconforming method each element interacts only with elements that share a face, 
whereas in the conforming method interactions also occur through shared edges and vertices. 
Moreover, from Figure~\ref{fig:spy}, 
we can compare the sparsity patterns of these two discretizations,
We present only the sparsity patterns corresponding to the \texttt{voro2} mesh.
For $k=1$, the sparsity patterns of the two methods are similar.
However, for $k=2$ we notice a saddle-point structure for the nonconforming case and a more complex, 
triple-level saddle-point structure for the conforming one, 
i.e., the global matrices have the following structure
\[
S_{c} = \begin{bmatrix}
A &E &F &G \\
E^\top  &B &H &J \\
F^\top  &H^\top &C &K \\
G^\top  &J^\top &K^\top &D \\
\end{bmatrix}\qquad\text{and}\qquad
S_{nc} = \begin{bmatrix}
A &C \\
C^\top &B
\end{bmatrix}\,,
\]
for the conforming and nonconforming case, respectively.
This structure is primarily due to the ordering of the degrees of freedom.
Consider the conforming case, in \texttt{vem++} the degrees of freedom are organized as follows:
first the vertex degrees of freedom, second the ones associated with edges, faces and polyhedra.
Since the non conforming case has only degrees of freedom associated with faces and polyhedra (\texttt{D1} and \texttt{D2}),
the linear system arising from this discretization has a saddle-point structure:
the first matrix entries are the degrees of freedom of the type \texttt{D1}
and the last part of the matrix entries are the degrees of freedom \texttt{D2}.

